\documentclass[11pt]{article}
\usepackage{latexsym,amsfonts,amsmath,epsfig,amssymb,mathrsfs,enumerate}
\usepackage{threeparttable}
\usepackage[pdftex, bookmarksnumbered, bookmarksopen, colorlinks, citecolor=blue, linkcolor=blue]{hyperref}
\usepackage{subfigure}
\usepackage[ruled]{algorithm2e}
\usepackage{bm}
\usepackage{booktabs}
\usepackage{amsmath}
\usepackage{multirow}
\usepackage{epstopdf}
\usepackage{threeparttable}
\usepackage{tabularx}
\usepackage{longtable, booktabs}
\usepackage{cases}
\usepackage{graphicx}
\usepackage{amsfonts}
\usepackage{amssymb}
\usepackage{color}
\usepackage{multirow}
\usepackage{multicol}
\usepackage{url}
\usepackage{amsfonts}
\usepackage{amsfonts,amssymb}

\usepackage{mathrsfs}
\numberwithin{equation}{section} \numberwithin{figure}{section}
\numberwithin{table}{section}

\newtheorem{lemma}{Lemma}[section]
\newcounter{nextauthor}
\newtheorem{Proposition}{Proposition}[section]
\newcommand{\new}{\newcommand*}
\new{\newt}{\newtheorem}

\newt{corollary}{Corollary}[section]
\newt{theorem}{Theorem}[section]
\newt{assumption}{Assumption}[section]
\newt{example}{Example}[section]
\newt{remark}{Remark}[section]

\def\mathrm{\mbox}

\makeatletter
\def\@biblabel#1{#1}

\begin{document}

 \begin{center}
  \noindent 	\centerline{\Large\bf Two nonmonotone multiobjective  memory gradient algorithms  \noindent \footnote{ This research is supported by the National Natural Science Foundation of China(12271071, 11991024), the Basic and Advanced Research Project of Chongqing
(cstc2021jcyj-msxmX0300), the Team Project of Innovation Leading Talent
in Chongqing (CQYC20210309536), the Chongqing Talent Plan contract system project
(cstc2022ycjh-bgzxm0147) and the Chongqing University Innovation Research
Group Project (CXQT20014). J. C. Yao was supported by the Grant MOST 111-2115-
M-039-001-MY2.}}\\

\vskip 0.4in
Jian-Wen Peng\footnote{ School of Mathematical Sciences, Chongqing Normal University, Chongqing 401331, China. e-mail: jwpeng168@hotmail.com}, Jie-Wen Zhang\footnote{ School of Mathematical Sciences, Chongqing Normal University, Chongqing 401331, China.} and Jen-Chih Yao\footnote{Corresponding author.  Research Center for Interneural Computing, China  Medical University,
 Taichung 40402, Taiwan and Department of Applied Mathematics, National Sun Yat-sen University, Kaohsiung 804,
 Taiwan. Email: yaojc@mail.cmu.edu.tw}

\end{center}

\vskip 0.1in

\begin{minipage}{5.5in}
{\bf Abstract}
In this paper,  two types of nonmonotone  memory gradient algorithm    for solving unconstrained multiobjective optimization problems are introduced.  Under some suitable conditions, we show the convergence of the full sequence generated by the proposed algorithms to a weak Pareto optimal point. 
\end{minipage}

\baselineskip 0.29in \vskip 2.0pt

\section{Introduction}

The multiobjective optimization problem (in short, MOP)  contains several objective functions have to be minimized or maximized simultaneously. The (MOP)  is a very important optimization model because of  its  usefule  applications in wide areas such as engineering, finance and management science. Since variables that minimize or maximize all objective functions at once do not necessarily exist, we use the concept of Pareto optimality. \\

There exist various iterative techniques  for solving the multiobjective optimization problem, such as the Newton method \cite{10,11}, quasi-Newton methods \cite{12,13,14}, the steepest descent method \cite{15}, and conjugate gradient methods \cite{2,3,4}, etc. Among these, the Newton method and quasi-Newton methods are renowned for their effectiveness in solving (MOP) due to their rapid convergence rates. However, these   Newton  type methods  need to compute the  Hession   matrices of the objective functions or their approximate  matrices, making them less suitable for large-scale problems. To address this limitation, the limited memory BFGS method for solving (MOP) has been employed (see \cite{1}). In contrast, the steepest descent method solving (MOP) doesn't rely on  any matrix operations but tends to have a slower convergence rate. Consequently, many efforts have been directed towards accelerating the steepest descent method  for solving (MOP) in recent years, leading to the emergence of accelerated gradient methods. Accelerated gradient methods  for solving (MOP) encompass various approaches, which include the conjugate gradient methods  and the memory gradient methods.

In iterative methods, only the current iterative information is used to generate a new iterate at each iteration and the previous iterative information is ignored. This is a waste of information in the design of efficient algorithms. The conjugate Gradient (in short, CG) methods for solving (MOP) use previous one-step iterative information to generate the next iterate at each iteration.  In general, the search direction can be defined by
\begin{eqnarray}
d_{k}=\left\{\begin{array}{ll}
v(x_{k}), & \text{if}\ k=0, \\
v(x_{k})+\beta_{k}d_{k-1}, & \text {if}\ k\geq1,
\end{array}\right.\label{cg}
\end{eqnarray}

where $\beta_{k}$ is a parameter that determines some different conjugate gradient methods. In \cite{2}, $\beta_{k}$'s were regarded as the extended forms of the five classical scalar optimization parameters. Additionally, the multiobjective extensions of the nonlinear conjugate gradient methods   were also proposed by  Liu-Storey \cite{3} and and Hager-Zhang \cite{4}, respectively.
The benefit of CG methods   is that they can handle complex optimization problems without the computation and storage of some matrices required by Newton-type methods. As generalizations of the  CG  methods, memory gradient (in short, MG) methods share the same characteristics as the  CG methods. The primary distinction between the   MG  methods and the  CG  methods  is that the  MG  methods may effectively leverage the prior multi-step iterative information to build a new iterate at each iteration, making it useful to design new algorithms with quick convergence rates.

Recently, Chen et al. \cite{5} proposed a multiobjective memory gradient algorithm (in short, MMG) for solving MOP. This method always satisfies
 the sufficient descent condition and converges globally if the   Wolf  line search conditions are satisfied. The
 main work of \cite{5} is to extend the memory gradient algorithm from  scalar optimization case   in \cite{18} to multiobjective
 optimization.

Two different types of  nonmonotonic memory gradient algorithms for scalar optimization are introduced and researched in \cite{16} and \cite{17}. In the field of multiobjective
optimization, nonmonotone line search techniques are frequently utilized and have shown their superior effectiveness when compared to
monotone algorithms (see \cite{6,7,8}).

 The main work of this paper is to extend the non-monotonic memory
gradient algorithms given in \cite{16} and \cite{17} to the two types of nonmonotone multiobjective memory
gradient algorithms and establish their global convergence.


The paper is organized as follows. Section 2 presents some notions, definitions and preliminary results. In Section 3, we propose two types of nonmonotone   multiobjective  memory gradient   algorithms. In Section 5, we establish the global convergence of the proposed NMMG motheds.  Finally, in Section 5 we give some conclusions.

\section{Preliminaries}
\setcounter{equation}{0}

Throughtout this paper,  $\mathbb{R}$, $\mathbb{R}_{+}$, $\mathbb{R}_{++}$ and $\mathbb{Z_{+}}$ stand for the sets of real numbers, non-negative real numbers, strictly positive real numbers and   the set of positive integers, respectively. For $m \in\mathbb{Z_{+}}$,  we use  $\Xi_m$ to denote the set $\{1,2,\ldots,m\}$. Let $e=(1,1,\ldots,1)^{\top}\in\mathbb{R}^{m}$,  $\langle\cdot,\cdot\rangle$ stand for the inner product in $\mathbb{R}^{n}$ and $\|\cdot\|$ denote the norm which means that $\| x \|=\sqrt{\langle x,x\rangle}$ for $x\in\mathbb{R}^{n}$. For two vectors $u$ and $v$ in $\mathbb{R}^{n}$,   $u\leq v$ denotes that $u_{i}\leq v_{i}$ for all $i=\Xi_m$ and $u<v$, denotes that $u_{i}<v_{i}$ for all $i=\Xi_m$.
The norm of a real matrix $A=(A_{i,j})_{m\times n}\in \mathbb{R}^{n}$ is defined as
 \begin{eqnarray}
\|A\|=\max\limits_{x\neq 0}\frac{\|Ax\|_{\infty}}{\|x\|}=\max\limits_{i\in\Xi_m}\|A_{i,\cdot}\|=\max\limits_{i\in\Xi_m}(\sum \limits_{j=1}^{n}A^{2}_{i,j})^{\frac{1}{2}}.\label{3}
 \end{eqnarray}

Let $a^{+} (a\in \mathbb{R})$ be denoted as

$$
a^{+}=\left\{\begin{array}{l}
0, \text { if } a=0, \\
\frac{1}{a}, \text { otherwise. }
\end{array}\right.
$$
Clearly, $aa^{+}\leq 1$ and $aa^{+}=1$  only when $a\neq 0$.

In this paper, we are concerned with the following (MOP):
 \begin{eqnarray}
\min\limits_{x\in\mathbb{R}^{n}}F(x)=(F_{1}(x),F_{2}(x),\ldots,F_{m}(x))^{\top},\label{mop}
 \end{eqnarray}
where $F_{i}:\mathbb{R}^{n}\rightarrow \mathbb{R},i\in\Xi_m$, are continuously differentiable and the superscript $" \top "$denotes the transpose. Given $x=(x_{1},x_{2},\ldots,x_{n})\in\mathbb{R}^{n}$, the Jacobian of $F$ at $x$ is fined by
$$JF(x)=[\nabla F_{1}(x),\nabla F_{2}(x),\ldots,\nabla F_{m}(x)]^{\top}.$$
The image of $JF(x)$ is denoted as
$$Im(JF(x))=\{JF(x)d:d\in \mathbb{R}^{n}\}.$$
\textbf{Definition 2.1} \cite{9} A vector $x^{*}\in\mathbb{R}^{n}$ is called Pareto optimal for (MOP), if there exists no $x\in \mathbb{R}^{n}$ such that $F(x)\leq F(x^{*})$ and $F(x)\neq F(x^{*})$.\\
\textbf{Definition 2.2} \cite{9} A vector $x^{*}\in\mathbb{R}^{n}$ is called Pareto optimal for (MOP), if
$$Im(JF(x^{*}))\cap(-\mathbb{R}^{m}_{++})=\emptyset,$$
where $JF(x^{*})$ denotes the image set of the Jacobian of $F$ at $x^{*}$.\\
\textbf{Definition 2.3} \cite{9} A vector $d\in\mathbb{R}^{n}$ is called descent direction for $F$ at $x$, if
$$JF(x)d\in-\mathbb{R}_{++}^{m}.$$

Now, we define $\psi:\mathbb{R}^{n}\times \mathbb{R}^{n}\rightarrow \mathbb{R}$ as follows:
\begin{eqnarray}
\psi(x,d)=\max\limits_{i\in\Xi_m}\langle\nabla F_{i}(x),d\rangle.\label{5}
 \end{eqnarray}

From the previous discussion, we known that $\psi$ can express Pareto critical and descent direction, i.e.,

(i) $d\in\mathbb{R}^{n}$ is a descent direction for $F$ at $x\in\mathbb{R}^{n}$ if $\psi(x,d)<0$,

(ii) $x\in\mathbb{R}^{n}$  is Pareto critical if $\psi(x,d)\geq 0$ for any $d\in\mathbb{R}^{n}$.

The following proposition illustrates several useful results related to $\psi$.
\begin{Proposition} \cite{5}For all $x,y\in\mathbb{R}^{n},\varrho>0$ and $b_{1},b_{2}\in\mathbb{R}^{n}$, we obain\\
 (i)$\psi(x,\varrho b_{1})=\varrho \psi(x,b_{1})$;\\
 (ii)$\psi(x,b_{1}+b_{2})\leq\psi(x,b_{1})+\psi(x,b_{2})$;\\
 (iii)$|\psi(x,b_{1})-\psi(x,b_{2})|\leq\|JF(x)b_{1}-JF(x)b_{2}\|$.
\end{Proposition}

Let us now consider the following scalar optimization problem:
 \begin{eqnarray}
\min\limits_{d\in\mathbb{R}^{n}}\psi(x,d)+\frac{1}{2}\|d\|^{2}. \label{6}
 \end{eqnarray}
Obviously, the objective of (\ref{6}) is   strongly convex. We denote the optimal solution and optimal value of (\ref{6}) by $v(x)$ and $\theta(x)$, respectively. Therefore
 \begin{eqnarray}
v(x)=\mathop{\arg\min}\limits_{d\in \mathbb{R}^{n}}\psi(x,d)+\frac{1}{2}\|d\|^{2},\label{7}
 \end{eqnarray}
and
 \begin{eqnarray}
\theta(x)=\psi(x,v(x))+\frac{1}{2}\|v(x)\|^{2}.\label{8}
 \end{eqnarray}
Indeed, problem (\ref{6}) can be equivalently rewritten as the following smooth quadratic programming problem:

\begin{equation}\begin{aligned}\min_{(t,d)\in\mathbb{R}\times\mathbb{R}^n}t+\frac{1}{2}\|d\|^2,\\\text{s.t. }\langle\nabla F_i(x),d\rangle\leq t,\ &i=\Xi_m. \label{09}\end{aligned}\end{equation}

Note that (\ref{09}) is a convex optimization problem with  linear constraints, then the strong dual holds. The Lagrangian function of (\ref{09}) is
$$L((t,d),\lambda)=t+\frac{1}{2}\|d\|^2+\sum_{i=1}^{m}\lambda_i(\langle\nabla F_i(x),d\rangle-t).$$
Under the Karush-Kuhn-Tucker (KKT) condition, we have
 \begin{eqnarray}
\sum\limits_{i =1}^{m}\lambda_{i}(x)=1,\label{ttt}
 \end{eqnarray}
 \begin{eqnarray}
v(x)+\sum_{i=1}^{m}\lambda_i(x)\nabla F_i(x)=0,\label{d0}
 \end{eqnarray}
 \begin{eqnarray}
\langle\nabla F_i(x),v(x)\rangle\le \theta(x),\ i\in \Xi_m,\label{lam0}
 \end{eqnarray}
\begin{eqnarray}
\lambda_i(x)\geq0,\  i\in \Xi_m,\label{laam0}
 \end{eqnarray}
\begin{eqnarray}
\lambda_i(\langle\nabla F_i(x),v(x)\rangle-\theta(x))=0,\  i\in \Xi_m,\label{scfb}
 \end{eqnarray}
From (\ref{d0}), we can obtain
 \begin{eqnarray}
v(x)=-\sum\limits_{i=1}^{m}(\lambda_i(x))\nabla F_{i}(x).\label{vxx}
 \end{eqnarray}
where $\lambda(x)=(\lambda_{1}(x),\lambda_{2}(x),\ldots,\lambda_{m}(x))$ is the solution to the dual problem:
\begin{equation}\begin{aligned}-\min_\lambda\frac12\|\sum_{i=1}^{m}\lambda_i\nabla F_i(x)\|^2,\\\text{s.t. }\lambda\in\Delta_m. \label{lamx}\end{aligned}\end{equation}
Recall that the strong dual holds and we get
 \begin{eqnarray}
\theta(x)=-\frac{1}{2}\|\sum_{i=1}^{m}\lambda_i(x)\nabla F_i(x)\|^2=-\frac{1}{2}\|v(x)\|^2.\label{cta}
 \end{eqnarray}

Let us now give a characterization of Pareto critical points of problem (\ref{mop}), which will be used in our
subsequent analysis.
\begin{Proposition} \cite{5} Let $v(\cdot)$ and $\theta(x)$ be defined as in (\ref{7}) and (\ref{8}), respectively. The following statements hold:\\
(i) If $x$ is a Pareto critical point of problem (\ref{mop}), then $v(x) = 0$ and 
$\theta(x)=0$; \\
(ii) If x is not a Pareto critical point of problem (\ref{mop}), then $v(x)\neq0,\theta(x)<0$ and $\psi(x,v(x))<-\|v(x)\|^{2}/2<0$;\\
(iii) $v(\cdot)$is continuous.
\end{Proposition}

\section{Nonmonotone multiobjective memory gradient algorithm}
This section describes two types of non-monotone multiobjective memory gradient algorithms (in short, NMMG)  for solving the problem (\ref{mop}). First, we give the (NMMG) algorithm with max-type Armijo line searches as follows:\\

\textbf{Algorithm 3.1:  Max-type NMMG algorithm}\\

Step 0.  Choose initial points $x_{0}\in\mathbb{R}^{n},\gamma_{0}>0$ and $N\in\mathbb{Z}_{+}$, parameters $0<\lambda_{1}\leq\lambda_{2},$ $0<\lambda_{3}\leq\lambda_{4}\leq1,\ \rho\in(0,1)$, and a nonnegative integer $M$.  Set $k:=0$.\\

Step 1.  Compute $v(x_{k})$ by using (2.13) and (2.14).

Step 2. If $v(x_{k})=0$ then STOP. Otherwise, continue with step 3.\\

Step 3.  Compute the search direction $d_{k}$ by the following definition of $d_{k}$
\begin{eqnarray}
d_{k}=\left\{\begin{array}{ll}
\gamma_{k}v(x_{k}), & \text{if}\ k=0, \\
\gamma_{k}v(x_{k})+\sum\limits_{j=1}^{N_{k}}\beta_{kj}d_{k-1}, & \text {if}\ k\geq1,
\end{array}\right.\label{dd}
\end{eqnarray}
and
\begin{eqnarray}
\beta_{kj}=-\frac{1}{N_{k}}\psi(x_{k},v(x_{k}))\phi^{+}_{kj},\label{12}
 \end{eqnarray}
where $\gamma_{k}>0$ and $\beta_{kj}\in\mathbb{R}(j\in\langle N_{k} \rangle,N_{k}=\min\{k,N\})$ are  algorithmic parameters.\\

Step 4. Compute a stepsize $\alpha_{k}$, let $M(k)=\min(k, M),\alpha_{k}^{0}\in[\lambda_{1},\lambda_{2}]$ and set $i:=0$.  if
 \begin{eqnarray}
F_{i}(x_{k}+\alpha_{k}d_{k})\leq \max\limits_{0\leq j \leq M(k)}F_{i}(x_{k-j})+\rho\alpha_{k}^{i}\psi(x_{k},d_{k}).\label{maxa1}
 \end{eqnarray}
holds, set $\alpha_{k}\equiv \alpha_{k}^{i}$. Otherwise go to step 5.\\

Step 5. Choose $\sigma_{k}^{i}\in[\lambda_{3},\lambda_{4}]$ and compute $\alpha_{k}^{i}$ such that
\begin{eqnarray}
\alpha_{k}^{i+1}=\alpha_{k}^{i}\sigma_{k}^{i}.\label{sigma}
 \end{eqnarray}\\

Step 6. Set $i:=i+1$ and go to step 4\\

Step 7. Let $x_{k+1}=x_{k}+\alpha_{k}d_{k}$, $k:=k+1$ and return to Step 1.\\

{\bf Remark 3.1} In Algorithm 3.1, we give a multiobjective version of a non-monotone memorized gradient algorithm that employs an max type Armijo line search technique. In single-objective optimization problems, this line search method speeds up convergence, and we expect the same effect for multi-objective optimization.\\

In scalar optimization, a new nonmonotonic technique called average type Armijo line search is given in \cite{8}. Therefore we will extend this non-monotonic technique in
\cite{8} to give the following non-monotone multiobjective    memory gradient algorithm with average type Armijo line searches:\\

\textbf{Algorithm 3.2:    Average-type NMMG algorithm}\\

Step 0.  Choose initial points $x_{0}\in\mathbb{R}^{n},\gamma_{0}>0$ and $N\in\mathbb{Z}_{+}$, parameters $\delta\in(0,1),\rho\in(0,1),$ $\tau_{k}=-\psi(x_{k},d_{k})/\|d_{k}\|^{2}, 0\leq_{\min}\leq\eta_{\max}\leq1$, $C_{0}=F(x_{0})$ and $Q_{0}=1$. Set $k:=0$.

Step 1. Compute $v(x_{k})$ by using (2.13) and (2.14).\\

Step 2. If $v(x_{k})=0$ then STOP. Otherwise, continue with step 3.\\

Step 3. Compute a search direction $d_{k}$ by (\ref{dd}) and (3.2).\\

Step 4. Compute the step length  $\alpha_{k}$:  Let $h_{k}$ is the smallest nonnegative integer $j$ such that
 \begin{eqnarray}
F_{i}(x_{k}+\tau_{k}\delta^{j} d_{k})\leq C_{k}^{i}+\rho\tau_{k}\delta^{j}\psi(x_{k},d_{k}).\label{maxa}
 \end{eqnarray}
 Let $\alpha_{k} =   \delta^{h_{k}}\tau_{k}$.
 \\

Step 5. Set $x_{k+1}=x_{k}+\alpha_{k}d_{k}$ .\\

Step 6. Choose $\eta_{k}\in[\eta_{\min},\eta_{\max}]$, update $Q_{k}$ and $C_{k}$ as follows:
\begin{eqnarray}
Q_{k+1}=\eta_{k} Q_{k}+1,\label{qk}
 \end{eqnarray}
\begin{eqnarray}
C_{k+1}^{i}=\frac{\eta_{k} Q_{k}+F_{i}(x_{k+1})}{Q_{k+1}}.\label{ave}
 \end{eqnarray}\\

Step 7. Set $k:=k+1$ and return to Step 1.

{\bf Remark 3.2}  It can be easily seen that $C_{k}^{i}$ is a convex combination of $F_{i}(x_{0}),F_{i}(x_{1}),\ldots,F_{i}(x_{k})$ for $i=\Xi_m$. Also, the choice of $\eta_{k}$ controls the degree of non-monotonicity of Algorithm 3.2. If $\eta_{k}=0$, for all $k$, then $C^{i}_{k} = F_{i}(x_{k})$ and the line search degenerates into the usual monotone Armijo line search. If $\eta_{k}=1$, for all $k$, then $C_{k} = A_{k}$, where
$$A_{k}=\frac{1}{k+1}\sum\limits_{i=0}^{k}F(x_{i})$$
is a vector of average function values. Moreover, as can be known in the next lemma, for each $i$, the values $F_{i} (x_{k})$ and $A_{k}^i$  are the lower and upper bounds of $C_{k}^{i}$, respectively, which also implies that the line search process is well-defined.

{\bf Remark 3.3}   If $N=1$ and $\gamma_{k}=1$ for every $k$, then   both  Algorithms 3.1 and 3.2  reduce to multiobjective nonlinear conjugate gradient algorithms in \cite{1}. If $N=1$ and $\gamma_{k}=1$ for each $k$ and $\beta_{kj}=0$ for all $k\geq 1$, then (\ref{dd}) can be viewed as the iterative form of the multiobjective steepest descent method in  \cite{15}.

\begin{lemma} Let $\{x_{k}\}$ be the sequences generated by  Algorithm 3.2. Then, we have
$$F(x_{k})\leq C_{k}\leq A_{k}.$$
Moreover, if $x_{k}$ is a noncritical point for (MOP), then there exists $\alpha_{k}$ satisfying the nonmonotone Armijo conditions of the line search procedure.
\end{lemma}

{\bf Proof.} Defining $D^{k}:\mathbb{R}\rightarrow \mathbb{R}^{m}$ as
$$D^{k}(t)=\frac{tC_{k-1}+F(x_{k})}{t+1}.$$
Then, its derivative is given by
$$(D^{k})^{'}(t)=\frac{\partial D^{k}(t)}{\partial t}=\frac{(C_{k-1}-F(x_{k}))}{(t+1)^{2}}$$
Since $JF(x_{k})d_{k}\leq 0$, it follows from the Armijo-type condition (\ref{maxa}), that $F(x_{k})\leq C_{k-1}$. This implies that $(D^{k})^{'}(t)\geq 0$ for all $t\neq -1$, that is, $D^{k}$ is nondecreasing for all $t\geq 0$. Also, from (\ref{qk}) and the fact that $\eta\in[0,1]$ and $Q_{0}=1$,we have $Q_{k}\geq 1$ for all $k$. Hence, $\eta Q_{k-1}\geq0$ holds, and we obtain
$$F(x_{k})=D^{k}(0)\leq D^{k}(\eta Q_{k-1})=C_{k}.$$

Now, let us prove $C_{k}\leq A_{k}$ by induction. For $k=0$, this inequality holds because $C_{0}=A_{0}=F(x_{0})$. So, assume that $C_{j}\leq A_{j}$ for all $0\leq j \leq k$. Because $\eta\in[0,1]$, $Q_{0}=1$  and (\ref{qk}), we can obtain
$$Q_{k+1}=1+\sum\limits_{i=0}^{k}\prod\limits_{l=0}^{i}\eta_{k-l}\leq k+2,$$
for each $k$. Thus, $0\leq Q_{k}-1\leq k$ holds. Since $D^{k}$ is nondecreasing for all $t\geq 0$ and $Q_{k}=\eta Q_{k-1}+1$ in Step 6, we obtain
$$C_{k}=D^{k}(\eta Q_{k-1})=D^{k}(Q_{k}-1)\leq D^{k}(k).$$
By the induction step, we also have
$$D^{k}(k)=\frac{kC_{k-1}+F(x_{k})}{k+1}\leq \frac{kA_{k-1}+F(x_{k})}{k+1}=A_{k},$$
and the conclusion follows.\\

As in the method given in \cite{5}, we denfine $\beta_{kj}(j\in \Xi_{ N_{k}})$ as follows
 \begin{eqnarray}
\beta_{kj}=-\frac{1}{N_{k}}\psi(x_{k},v(x_{k}))\phi^{+}_{kj},\label{12}
 \end{eqnarray}
where $\phi_{kj}(j\in \Xi_{ N_{k}})$ are parameters satisfying the following relation:
\begin{eqnarray}
\phi_{kj}>\max\{\frac{\psi(x_{k},d_{k-j})}{\gamma_{k}},0\}.\label{13}
 \end{eqnarray}
Note that $\beta_{kj}>0$ since $\psi(x_{k},v(x_{k}))<0$ for any $k$.

The subsequent property shows that if the relevant parameters satisfy (\ref{12}) and (\ref{13}), then $d_{k}$ is a descending direction.
\begin{lemma}  \cite{5}Let the direction $d_{k}$ be given in (\ref{dd}). Assume that $\beta_{kj}$ and $\phi_{kj}$ satisfy (\ref{12}) and (\ref{13}) for $k\geq1$ and $j\in\Xi_{ N_{k}}$, respectively. Then, $d_{k}$ is a descent direction for all $k$.
\end{lemma}

The following conditions must be mentioned in the convergence analysis of scalar optimization,
 \begin{eqnarray*}
\langle\nabla F_{1}(x_{k}),d_{k}\rangle\leq-c\left\|\nabla F_{1}(x_{k})\right\|^{2},
 \end{eqnarray*}
which is the well-known sufficient descent condition. Likewise, in convergence analysis for (MOP), we will need the more stringent condition,
\begin{eqnarray}
\psi(x_{k},d_{k})\leq c\psi(x_{k},v(x_{k})).\label{17}
 \end{eqnarray}
If (\ref{17}) holds, we say that a direction $d_{k}$ satisfies the sufficient descent condition at $x_k$. It is worth noting that the
general concept of sufficient descent condition was first introduced by Lucambio P\'{e}rez and Prudente \cite{2} in the field of
vector optimization. We provide a sufficient condition on $\gamma_{k}$ and $\phi_{kj}$ in the following lemma to guarantee the descent property on $d_{k}$.

\begin{lemma}
\cite{5}  Let the direction $d_{k}$ be given in (\ref{dd}). Suppose that there is a positive constant $\gamma^{*}$ such that $\gamma_{k}\geq\gamma^{*},\beta_{kj}$ satisfies (\ref{13}) and $\phi_{kj}$ has the following property:
 \begin{eqnarray}
\phi_{kj}>\frac{\psi(x_{k},d_{k-j})+\left\|JF(x_{k})\right\|\left\|d_{k-j}\right\|}{\gamma_{k}}.\label{18}
 \end{eqnarray}
 Then, $d_{k}$ satisfies the sufficient descent condition (\ref{17}) with $\sigma=\gamma^{*}/2>0$ for any $k$.
\end{lemma}


In the next lemma, we will show that both Algorithms 3.1 and 3.2 are well-defined.\\

\begin{lemma}  Let $\{x_{k}\}$ be an iterate of Algorithms 3.1 ( or   3.2, respectively). If $x_{k}$ is not Pareto critical, then there exists a stepsize $\alpha_{k}>0$ satisfying the Armijo-type conditions (\ref{maxa1})(or (\ref{maxa}), respectively).
\end{lemma}
{\bf Proof.} Since $x_{k}$ is not a Pareto critical point, $\psi(x_k, d_k)<0$  holds. Furthermore, $\psi(x_k, d_k)<\delta\psi(x_k, d_k)$ because $\delta\in(0,1)$. Now, since $F$ is differentiable, we also have

 \begin{eqnarray}
F(x_{k}+\alpha d_{k})=F(x_{k})+\alpha JF(x_{k})d_{k}+o(\alpha).\label{23}
 \end{eqnarray}
 By (\ref{23}), we obtain
\begin{equation}
	\begin{split}
	F(x_{k}+\alpha d_{k})&= F(x_{k})+\alpha JF(x_{k})d_{k}+o(\alpha)\\
	        &\leq F(x_{k})+\alpha \delta JF(x_{k})d_{k}\\
&\leq F(x_{k})+\alpha\delta\psi(x_{k},d_{k}),
	\end{split}\label{23.1}
\end{equation}
Therefore, there exists $\bar{\alpha}\in(0, 1)$ such that
$$F(x_{k}+\alpha d_{k})\leq F(x_{k})+\alpha \delta\psi(x_{k},d_{k}),\forall \alpha\in(0,\bar{\alpha}].$$
Since $F(x_{k} )\leq C_{k}$, the above inequality satisfies
$$F(x_{k}+\alpha d_{k})\leq C_{k}+\alpha \delta\psi(x_{k},d_{k}), \forall \alpha\in(0,\bar{\alpha}],$$
or
$$F(x_{k}+\alpha_{k}d_{k})\leq \max\limits_{0\leq j \leq M(k)}F(x_{k-j})+\rho\alpha_{k}\psi(x_{k},d_{k}), \forall\alpha\in(0,\bar{\alpha}],$$
and the proof is finished.

\section{Global convergence   of  (NMMG) algorithms} \setcounter{equation}{0}
In this section, we show the global convergence property of the two types of (NMMG) algorithms. For this purpose, we need the following assumptions:\\
{\bf (A1)} $F$ is bounded from below on the level set $\mathcal{L}=\{x\in\mathbb{R}^{n}:F(x) \leq F(x_{0})\}$, where $x_{0}\in\mathbb{R}^{n}$ is an available  point.\\
{\bf (A2)} The Jacobian $JF$ is Lipschitz continuous in an open convex set $\mathcal{B}$ that contains $\mathcal{L}$, i.e., there exists a constant $L>0$ such that $\left\|JF(x)-JF(y)\right\|\leq L\left\|x-y\right\|$ for all $x,y\in\mathcal{B}$.\\
{\bf (A3)} The set $\mathcal{L}$ is bounded.\\
{\bf Remark 4.1} Under the assumption   (A3), the sequence $\{ v(x_{k}),d_{k} \}$ is bounded. By \cite{5}, it is known that there exist constants $\xi_{1},\xi_{2}>0$ such that $\|v(x_{k})\|\leq\xi_{1}$ and $\|d_{k}\|\leq\xi_{2}$. Therefore, $\langle v(x_{k}),d_{k}\rangle\leq \|v(x_{k})\| \|d_{k}\| \leq \xi_{1}\xi_{2}$.\\


Now we establish  the global convergence of Algorithm 3.1.

\begin{theorem}
Assume that the sequence $\{(x_{k}, d_{k})\}$ is produced by Algorithm 3.1. If the assumptions  (A1), (A2) and (A3) hold true, then $\lim\inf_{k\rightarrow\infty}\|v(x_{k})\|=0$.
\end{theorem}
{\bf Proof.} Define $l(k)$ as a number such that
$$k-M(k)\leq l(k)\leq k\ and\ F^{l(k)}=\max\limits_{0\leq  j \leq M(k)}\{F^{k-j}\},$$
for any $k$. By (\ref{maxa1}) and the fact that $M(k+1)\leq M(k)+1$, we have
\begin{equation*}
	\begin{split}
	F^{l(k+1)}&=\max\limits_{0\leq  j \leq M(k+1)}\{F^{k+1-j}\}\\
	        &\leq \max\limits_{0\leq  j \leq M(k)+1}\{F^{k+1-j}\} \\
&=\max\{ F^{l(k)},F^{k+1}\}\\
&=F^{l(k)},
	\end{split}
\end{equation*}
Thus, we see that the sequence $\{  F^{l(k)} \}$ does not increase. Moreover, by (\ref{maxa1}), we get that
\begin{equation*}
	\begin{split}
	F^{l(k+1)}&\leq F^{l(l(k))}\\
	        &\leq \max\limits_{0\leq  j \leq M(l(k)-1)}\{F^{l(k)-1-j}+\rho\alpha_{l(k)-1}\psi(x_{l(k)-1},d_{l(k)-1})\} \\
&=F^{l(l(k)-1)}+\rho\alpha_{l(k)-1}\psi(x_{l(k)-1},d_{l(k)-1})\}.
	\end{split}
\end{equation*}

By assumptions (A1) and (A2) and the fact that the sequence $\{  F^{l(k)} \}$ is nonincreasing, $\{  F^{l(k)} \}$ has a limit. In the remainder of the proof, we replace the subsequence $\{  l(k)-1 \}$ by $\{  k^{'} \}$. Thus, we have
\begin{equation}\label{al1}
  \lim\limits_{k^{'}\rightarrow \infty}\ \alpha_{k^{'}}\psi(x_{k^{'}},d_{k^{'}})=0.
\end{equation}

If the theorem does not hold, then there exists a constant $\gamma >0$ such that
 \begin{eqnarray}
\|v(x_{k})\|\geq \gamma\label{al2}
 \end{eqnarray}
for any $k$.

Since  $\theta(x_k)=\psi(x_k,v(x_k))+\frac{1}{2}\|v(x_k)\|^{2} \leq 0$, we have
 \begin{eqnarray}
 \alpha_{k^{'}}\psi(x_{k^{'}},d_{k^{'}})\leq  -\frac{1}{2}\alpha_{k^{'}}\|v(x_{k^{'}})\|^{2}\leq  -\frac{1}{2}\alpha_{k^{'}}\gamma^{2}<0.\label{al3}
 \end{eqnarray}
It follows From (\ref{al1}) and (\ref{al3})  that
$$ \lim\limits_{k^{'}\rightarrow \infty}\ \alpha_{k^{'}}=0.$$
This equation implies that when $k^{'}$ is very large, $\alpha_{k^{'}}^{i}$ does not satisfy (\ref{maxa1}), i.e. $\alpha_{k^{'}}^{i}$ holds for some $i\neq 0$, and $\alpha_{k^{'}}^{i-1}$ do not satisfy (\ref{maxa1}). Thus,   we get
\begin{equation}
	\begin{split}
	F(x_{k^{'}}+\alpha_{k^{'}}^{i-1}d_{k^{'}})&> \max\limits_{0\leq j M(k^{'})}\{ F^{k^{'}-j} \}+\rho\alpha_{k^{'}}^{i-1}\psi(x_{k^{'}},d_{k^{'}})\\
	        &\geq F^{k^{'}}+\rho\alpha_{k^{'}}^{i-1}\psi(x_{k^{'}},d_{k^{'}}),
	\end{split}\label{al4}
\end{equation}
that is
\begin{equation}\label{al5}
  F(x_{k^{'}}+\alpha_{k^{'}}^{i-1}d_{k^{'}})-F^{k^{'}}\geq \rho\alpha_{k^{'}}^{i-1}\psi(x_{k^{'}},d_{k^{'}}).
\end{equation}
Using the mean-value theorem and the Lipschitz continuity of JF, $t\in(0,1)$ we obtain that
\begin{equation}
	\begin{split}
	 F_{i}(x_{k^{'}}+\alpha_{k^{'}}^{i-1}d_{k^{'}})-F^{k^{'}}&=\alpha_{k^{'}}^{i-1}\nabla F_{i}(x_{k^{'}}+t\alpha_{k^{'}}^{i-1}d_{k^{'}})^{\top}d_{k^{'}} \\
	        &=\alpha_{k^{'}}^{i-1}[\nabla F_{i}(x_{k^{'}})^{\top}d_{k^{'}}+(\nabla F_{i}(x_{k^{'}}+t\alpha_{k^{'}}^{i-1}d_{k^{'}})-\nabla F_{i}(x_{k^{'}}))^{\top}d_{k^{'}}]\\
&\leq \alpha_{k^{'}}^{i-1}\{\nabla F_{i}(x_{k^{'}})^{\top}d_{k^{'}}+Lt\alpha_{k^{'}}^{i-1}\|d_{k^{'}}\|^{2}  \},
	\end{split}\label{al6}
\end{equation}
It follows from (\ref{sigma})and (\ref{al5})that
\begin{equation}\label{al7}
  \alpha_{k^{'}}^{i-1}\{\nabla F_{i}(x_{k^{'}})^{\top}d_{k^{'}}+Lt\alpha_{k^{'}}^{i-1}\|d_{k^{'}}\|^{2}  \} > \rho\alpha_{k^{'}}^{i-1}\psi(x_{k^{'}},d_{k^{'}}).
\end{equation}
Thus, we obtain
\begin{equation}\label{al8}
 \nabla F_{i}(x_{k^{'}})^{\top}d_{k^{'}}+Lt\alpha_{k^{'}}^{i-1}\|d_{k^{'}}\|^{2}   > \rho\psi(x_{k^{'}},d_{k^{'}}),
\end{equation}
Since $\psi(x_{k},d_{k})\geq \nabla F_{i}(x_{k})_{\top}d_{k}$, we can get
\begin{equation}\label{al19}
Lt\alpha_{k^{'}}^{i-1}\|d_{k^{'}}\|^{2} > \rho\psi(x_{k^{'}},d_{k^{'}})- \nabla F_{i}(x_{k^{'}})^{\top}d_{k^{'}} \geq (\rho-1)\psi(x_{k^{'}},d_{k^{'}}).
\end{equation}
It follows from $\alpha_{k^{'}}^{i+1}=\alpha_{k^{'}}^{i}\sigma_{k^{'}}^{i}$ and (\ref{al19}) that
\begin{equation*}
Lt\frac{\alpha_{k^{'}}}{\sigma_{k^{'}}^{i-1}}\|d_{k^{'}}\|^{2} > \rho\psi(x_{k^{'}},d_{k^{'}})- \nabla F_{i}(x_{k^{'}})^{\top}d_{k^{'}} \geq (\rho-1)\psi(x_{k^{'}},d_{k^{'}}).
\end{equation*}
Considering the conditions $\delta\in(0,1)$ and $\lambda_{3}\leq\sigma_{k^{'}}^{i-1}$, we can write
\begin{equation*}
\alpha_{k^{'}}> \frac{\sigma_{k^{'}}^{i-1}(\rho-1)\psi(x_{k^{'}},d_{k^{'}})}{Lt\|d_{k^{'}}\|^{2}}>\frac{\lambda_{3}(\rho-1)\psi(x_{k^{'}},d_{k^{'}})}{Lt\|d_{k^{'}}\|^{2}}.
\end{equation*}
As a result, we get
\begin{equation}\label{al20}
\alpha_{k^{'}}>\bar{c}\frac{|\psi(x_{k^{'}},d_{k^{'}})|}{\|d_{k^{'}}\|^{2}},
\end{equation}
where $\bar{c}=\frac{\lambda_{3}(1-\rho)}{Lt}$.

It follows from  (\ref{al20}) that
$$\alpha_{k^{'}}|\psi(x_{k^{'}},d_{k^{'}})|> \bar{c} \frac{\psi^{2}(x_{k^{'}},d_{k^{'}})}{\|d_{k^{'}}\|^{2}}.$$
From (\ref{al1}) and the above inequality, we have
\begin{equation}\label{al21}
  \lim\limits_{k^{'}\rightarrow \infty}\ \frac{\psi^{2}(x_{k^{'}},d_{k^{'}})}{\|d_{k^{'}}\|^{2}}=0.
\end{equation}
On the other hand, observing that $0<\gamma^{2}\leq\|v(x_{k})\|^{2}\leq-2\psi(x_{k},v(x_{k}))\leq -2\psi(x_{k},d_{k})/c$, where the last inequality follows from (\ref{17}) we get
$$|\psi(x_{k},d_{k})|\geq \frac{c\gamma^{2}}{2}.$$
Thus, according to Remark 3.4, the same formula as Eq. (37) in \cite{5} can be obtained, namely
$$\langle v(x_{k}),d_{k}\rangle\leq a |\psi(x_{k},d_{k})|,$$
where $a=2\xi_{1}\xi_{2}/(c\gamma^{2})$.

Therefore, by similar arguments with that in  Theorem 4.1 of \cite{5} leads to the  conclusion that

\begin{equation}\label{al22}
  \frac{\psi^{2}(x_{k^{'}},d_{k^{'}})}{\|d_{k^{'}}\|^{2}}\geq \frac{\|v(x_{k})\|^{2}}{4+a^{2}}\geq \frac{\gamma^{2}}{4+a^{2}}>0
\end{equation}
which contradicts (\ref{al21}). Thereforethe proof is complete.

\begin{lemma}   Let $\{x_{k}\}$ be the sequence generated by Algorithm 3.2. Then, $\{C_{k}^{i}\}$ is nonincreasing and admits a limit when $k\rightarrow\infty$.
\end{lemma}

{\bf Proof.} Since $d_{k}$ is a descent direction, $\psi(x_{k}),d_{k}<0$ for all $k$. Then, the Armijo-type condition (\ref{maxa}) shows that $F_{i}(x_{k+1})\leq C_{k}^{i}$ for all $i$. It is clear from (\ref{ave}) that
\begin{equation*}
	\begin{split}
	C_{k+1}^{i}&= \frac{\eta_{k} Q_{k}}{Q_{k+1}}C_{k}^{i}+\frac{1}{Q_{k+1}}F_{i}(x_{k+1})\\
	        &\leq \frac{\eta_{k} Q_{k}}{Q_{k+1}}C_{k}^{i}+\frac{C_{k}^{i}+\rho\tau_{k}\delta^{h_{k}}\psi(x_{k},d_{k})}{Q_{k+1}}\\
&=C_{k}^{i}+\frac{\rho\tau_{k}\delta^{h_{k}}\psi(x_{k},d_{k})}{Q_{k+1}}\\
&\leq C_{k}^{i},
	\end{split}
\end{equation*}
which also implies that $\{C_{k}^{i}\}$ is nonincreasing. Since $F_{i}$ is bounded from below, and $F_{i}(x_{k})\leq C_{k}^{i}$ for all $i$ and $k$, we can conclude that $\{C_{k}^{i}\}$ admits a limit when $k\rightarrow \infty$.\\

The following  lemmas display the decrease property for the function value of iterate points generated by Algorithm 3.2.\\

\begin{lemma} Assume that the sequence $\{(x_{k}, d_{k})\}$ is produced by Algorithm 3.2 and that (A1) and (A2) hold true. Then, there is a positive constant $\omega$ such that
 \begin{eqnarray}
C_{k}-F(x_{k+1})\geq \omega\frac{\psi^{2}(x_{k},d_{k})}{\|d_{k}\|^{2}}e,\ \forall k.\label{25}
 \end{eqnarray}
\end{lemma}
{\bf Proof.}  we have the following two cases.

Case 1. Let $\alpha_{k}$ satisfy the line search of Algorithm 3.2, we have
\begin{equation}
	\begin{split}
	F_{i}(x_{k+1})&\leq C_{k}^{i}+\rho\alpha_{k}\psi(x_{k},d_{k})\\
	        &=C_{k}^{i}+\rho\tau_{k}\delta^{h_{k}}\psi(x_{k},d_{k})\\
&=C_{k}^{i}-\rho\delta^{h_{k}}\frac{\psi^{2}(x_{k},d_{k})}{\|d_{k}\|^{2}},
	\end{split}\label{ca1}
\end{equation}
Therefore, it is concluded that
$$C_{k}-F(x_{k+1})\geq \omega\frac{\psi^{2}(x_{k},d_{k})}{\|d_{k}\|^{2}}e,$$
where $\omega=\rho\delta^{h_{k}}$.

Case 2. By the non-monotonic Armijo line search of Algorithm 3.2, when $\alpha_{k} =
\max\{\tau_{k},\delta\tau_{k},\delta^{2}\tau_{k},\ldots\}$ satisfies the line search (\ref{maxa}),
$\alpha=\alpha_{k}/\delta$ does not satisfy  (\ref{maxa}), and we have
\begin{equation}
	\begin{split}
	F_{i}(x_{k}+\alpha d_{k})&> C_{k}^{i}+\rho\alpha\psi(x_{k},d_{k})\\
	        &\geq F_{i}(x_{k})+\rho\alpha\psi(x_{k},d_{k}),
	\end{split}\label{ca2}
\end{equation}
that is say
\begin{equation}\label{ca22}
  F_{i}(x_{k}+\alpha d_{k})-F_{i}(x_{k})>\rho\alpha\psi(x_{k},d_{k}).
\end{equation}

Applying the mean-value theorem on the left side of (\ref{ca22}), there exists $t\in[0,1]$, such that
$$\alpha\nabla F_{i}(x_{k}+t\alpha d_{k})^{\top}d_{k}>\rho\alpha\psi(x_{k},d_{k}),$$
i.e.,
\begin{equation}\label{ca3}
  \nabla F_{i}(x_{k}+t\alpha d_{k})^{\top}d_{k}>\rho\psi(x_{k},d_{k}).
\end{equation}
 By the definition of $\psi(x,d)$, we get $\psi(x_{k},d_{k})>\nabla F_{i}(x_{k})^{\top}d_{k}$ . Combining this with the above inequality (\ref{ca22}), we have
\begin{equation}
	\begin{split}
	(\rho-1)\psi(x_{k},d_{k})&<\nabla F_{i}(x_{k}+t\alpha d_{k})^{\top}d_{k}-\psi(x_{k},d_{k}) \\
	        &\leq \nabla F_{i}(x_{k}+t\alpha d_{k})^{\top}d_{k}-\nabla F_{i}(x_{k})^{\top}d_{k}\\
&=\langle \nabla F_{i}(x_{k}+t\alpha d_{k}),d_{k}\rangle-\langle \nabla F_{i}(x_{k}),d_{k}\rangle.
	\end{split}\label{ca31}
\end{equation}
By the Cauchy-Schwarz inequality, (A2) and $t\in[0,1]$, we get
\begin{equation}
	\begin{split}
	\langle \nabla F_{i}(x_{k}+t\alpha d_{k}),d_{k}\rangle-\langle \nabla F_{i}(x_{k}),d_{k}\rangle&\leq\|\nabla F_{i}(x_{k}+t\alpha d_{k})-\nabla F_{i}(x_{k})\| \|d_{k}\| \\
	        &\leq Lt\alpha\|d_{k}\|^{2}\\
&\leq L\alpha\|d_{k}\|^{2}.
	\end{split}\label{ca32}
\end{equation}
Thus, by (\ref{ca31}) and (\ref{ca32}), we get
\begin{equation}\label{ca33}
  \alpha_{k}\geq \frac{\delta(\rho-1)}{L}\frac{\psi(x_{k},d_{k})}{\|d_{k}\|^{2}}.
\end{equation}
From (\ref{maxa}), we get
\begin{equation}\label{ca34}
 C_{k}^{i}- F_{i}(x_{k+1})\geq -\rho\alpha_{k}\psi(x_{k},d_{k})\geq \frac{\delta\rho(1-\rho)}{L}\frac{\psi^{2}(x_{k},d_{k})}{\|d_{k}\|^{2}}.
\end{equation}
where $\omega=\frac{\delta\rho(1-\rho)}{L}$. Therefore, if we set
$$\omega=\{ \rho\delta^{h_{k}}, \frac{\delta\rho(1-\rho)}{L} \},$$
then the result of (\ref{25}) is satisfied.

\begin{lemma} Assume that the sequence $(x_{k}, d_{k})$ is produced by Algorithm 3.2. If (A1) and (A2) are satisfied, then
 \begin{eqnarray}
\sum\limits_{k\geq 0}\frac{1}{Q_{k+1}}\frac{\psi_{2}(x_{k},d_{k})}{\|d_{k}\|^{2}}<+\infty.\label{27}
 \end{eqnarray}
\end{lemma}

{\bf Proof.} By (\ref{ave}) and (\ref{25}), we have
\begin{equation*}
	\begin{split}
	C_{k+1}&=\frac{\eta_{k}Q_{k}C_{k}+F(x_{k+1})}{Q_{k+1}}\\
	        &\leq \frac{\eta_{k}Q_{k}C_{k}+C_{k}-\omega\frac{\psi^{2}(x_{k},d_{k})}{\|d_{k}\|^{2}}}{Q_{k+1}}\\
&=C_{k}-\frac{\omega}{Q_{k+1}}\frac{\psi^{2}(x_{k},d_{k})}{\|d_{k}\|^{2}}.
	\end{split}
\end{equation*}
Thus,
$$\frac{\omega}{Q_{k+1}}\frac{\psi^{2}(x_{k},d_{k})}{\|d_{k}\|^{2}}\leq C_{k}-C_{k+1}.$$
Therefor, we obtain
\begin{eqnarray}
\sum\limits_{k\geq 0}\frac{1}{Q_{k+1}}\frac{\psi_{2}(x_{k},d_{k})}{\|d_{k}\|^{2}}<+\infty.\label{37}
 \end{eqnarray}

Now, we will prove the global convergence of Algorithm 3.2. \\

\begin{theorem}
Assume that the sequence $\{(x_{k}, d_{k})\}$ is produced by Algorithm 3.2. If the assumptions (A1),(A2) and (A3) hold true, then $\lim\inf_{k\rightarrow\infty}\|v(x_{k})\|=0$.
\end{theorem}
{\bf Proof.} Suppose by contradiction that there is $\varsigma>0$ such that
$$\|v(x_{k}\|)\geq\varsigma$$
for any $k$.

 By Thorem 4.1 in \cite{5}, we obatian
$$\sum\limits_{k=0}^{\infty}\frac{\psi^{2}(x_{k},d_{k})}{\|d_{k}\|^{2}}\geq \infty$$
Obviously, by using the relation $\eta_{\max}<1$ and (\ref{qk}), we have
$$Q_{k+1}=1+\sum\limits_{i=0}^{k}\prod\limits_{l=0}^{i}\eta_{k-l}\leq 1+\sum\limits_{i=0}^{k}\eta_{\max}^{i+1}\leq \sum\limits_{i=0}^{k}\eta_{\max}^{i}=\frac{1}{1-\eta_{\max}}.$$
Thus,
$$\frac{1}{Q_{k+1}}\geq 1-\eta_{\max}.$$
So, we can obtain
\begin{eqnarray*}
\sum\limits_{k\geq 0}\frac{1}{Q_{k+1}}\frac{\psi^{2}(x_{k},d_{k})}{\|d_{k}\|^{2}}\geq (1-\eta_{\max})\sum\limits_{k\geq 0}\frac{\psi^{2}(x_{k},d_{k})}{\|d_{k}\|^{2}}>\infty,
 \end{eqnarray*}
which contradicts with (\ref{37}). Thus,    the global convergence of the Algorithm 3.2 is proved.

\section{ Conclusion}  \setcounter{equation}{0}

In this work, we propose two non-monotonic memory gradient algorithms to solve multiobjective optimization problems. The novelty of
the approach of these algorithms lies in that the classical monotonic strategy  has been replaced by   the two types of nonmonotone linear search rule, which is called max type Armijo line search or average type Armijo line search. Under
mild assumptions, we establish that the sequences generated by the proposed two  non-monotone  multiobjective  memory gradient algorithms converges global to weak Pareto optimal
solutions.

We leave the numerical implementation of the proposed nonmonotone memory gradient algorithm and its application to solving some
 practical problems for future work.

\end{document}